\documentclass[jsl]{asl}
\title{Presburger sets and p-minimal fields.}
\author[Raf Cluckers]{Raf Cluckers$^\dag$}
\revauthor{Cluckers, Raf}
\thanks{$\dag$Research Assistant of the Fund for Scientific Research --
 Flanders (Belgium)(F.W.O.)} \revauthor{Cluckers, Raf}
\address{Department of Mathematics\\
Katholieke Universiteit Leuven\\ Celestijnenlaan 200B\\ B-3001 Leuven,
Belgium}
 \email{raf.cluckers@wis.kuleuven.ac.be}
 \urladdr{http://www.wis.kuleuven.ac.be/algebra/raf/}
\newtheorem{theorem}{Theorem}
\newtheorem{proposition}{Proposition}
\newtheorem{lemma}{Lemma}

\newtheorem{cor}{Corollary}
\theoremstyle{definition}
\newtheorem{definition}{Definition}
\newtheorem*{remark}{Remark}
\newtheorem*{remarks}{Remarks}

\newcommand{\Z}{\mathbb{Z}}

\newcommand{\Q}{\mathbb{Q}}
\newcommand{\Pm}{\mathcal{P}}
\newcommand{\Lm}{\mathcal{L}}
\newcommand{\SSS}{\mathcal{S}}

\newcommand{\Lp}{\mathcal{L}_{Pres}}
\newcommand{\Lmac}{\mathcal{L}_{Mac}}
\DeclareMathOperator{\sq}{\square}
\begin{document}
\begin{abstract}
We prove a cell decomposition theorem for Presburger sets and
introduce a dimension theory for $Z$-groups with the Presburger
structure. Using the cell decomposition theorem we obtain a full
classification of Presburger sets up to definable bijection. We
also exhibit a tight connection between the definable sets in an
arbitrary p-minimal field and Presburger sets in its value group.
We give a negative result about expansions of Presburger
structures and prove  uniform elimination of imaginaries for
Presburger structures within the Presburger language.
\end{abstract}
\maketitle
\section{Introduction}

At the ``Alg\`ebre, Logique et Cave Particuli\`ere'' meeting in
Lyon (1995), A.~Pillay posed the question of whether there exists
some dimension theory for $Z$-groups with the Presburger structure
which would give rise to a classification of all Presburger sets
up to definable bijection, possibly using other invariants as
well. In this paper we answer this question of Pillay: we classify
the Presburger sets up to definable bijection
(Thm.~\ref{classification}), using as only classifying invariant
the (logical) algebraic dimension. In order to prove this
classification, we first formulate a cell decomposition theorem
for Presburger groups (Thm.~\ref{cell decomp}) and a
rectilinearisation theorem for the definable sets
(Thm.~\ref{recti}). Also a rectilinearisation theorem depending on
parameters is proven (Thm.~\ref{param recti}).

Expansions of Presburger groups have recently been studied
intensively. One could say that on the one hand one looks for
(concrete) expansions which remain decidable and have bounded
complexity, and on the other hand different kinds of minimality
conditions (like coset-minimality, etc.) are used to characterize
general classes of expansions (see e.g. \cite{Wagner},
\cite{Point}). In section \ref{section Michaux} we examine
expansions of Presburger groups satisfying natural kinds of
minimality conditions.

In \cite{Haskell}, D.~Haskell and D.~Macpherson defined the notion
of p-minimal fields, as a $p$-adic counterpart of $o$-minimal
fields. A p-minimal field always is a $p$-adically closed field,
and its value group is a $Z$-group. Interactions between definable
sets in a given $p$-adically closed field and Presburger sets in
its value group have been studied in the context of $p$-adic
integration for several p-minimal structures (see \cite{Denef1},
\cite{suban}). In Theorem \ref{application}, we exhibit a close
connection between definable sets in arbitrary p-minimal fields
and Presburger sets in the corresponding value groups.

In the last section, we use the cell decomposition theorem in an
elementary way to obtain uniform elimination of imaginaries for
Z-groups
 without introducing extra sorts.
\\
\subsection{Terminology and notation.}
In this paper $G$ always denotes a $Z$-group, i.e.~a group which
is elementary equivalent to the integers $\Z$ in the Presburger
language $\Lp = \langle
+,\leq,\{\equiv(\mathrm{mod}{n})\}_{n>0},0,1 \rangle$ where
$\equiv(\mathrm{mod}{n})$ is the equivalence relation in two
variables modulo the integer $n>0$. We call $(G,\Lp)$ a Presburger
structure and we write $H$ for the nonnegative elements in $G$. By
a Presburger set, function, etc., we mean a $\Lp$-definable set,
function, etc., and by \emph{definable} we always mean definable
with parameters (otherwise we say $\emptyset$-definable,
$S$-definable, etc.). We call a set $X\subset G^m$ \emph{bounded}
if there is a tuple $z\in H^m$ such that $-z_i\leq x_i\leq z_i$
for each $x\in X$ and $i=1,\ldots,m$. For $k\leq m$ we write
$\pi_k:G^m\to G^k$ for the projection on the first $k$ coordinates
and for $X\subset G^{k+n}$ and $x\in\pi_k(X)$ we write $X_x$ for
the fiber $\{y\in G^n\mid (x,y)\in X\}$. We recall that the theory
$\mathrm{Th}(\Z,\Lp)$ has definable Skolem functions, quantifier
elimination in $\Lp$ and is decidable \cite{Presburger}.
\section{Cell Decomposition Theorem}
We prove a cell decomposition theorem for Presburger structures,
by first proving it in dimension 1 and subsequently using a
compactness argument. An elementary arithmetical proof can also be
given, using techniques like in the proof of Lemma 3.2 in
\cite{Denef3}, but our proof has the advantage that it goes
through in other contexts as well (see section \ref{section
Michaux} and \ref{section p}). As always, $G$ denotes a $Z$-group.
 \begin{definition}\label{linearity}
We call a function $f:X\subset G^m\to G$ \emph{linear} if there is
a constant $\gamma\in G$ and integers $a_i$,  $0\leq c_i < n_i$
for $i=1,\ldots,m$ such that $x_i-c_i\equiv 0 \pmod{n_i}$ and
\[
f(x)=\sum_{i=1}^m a_i(\frac{x_i-c_i}{n_i})+\gamma.
\]
for all $x=(x_1,\ldots,x_m)\in X$. We call $f$ \emph{piecewise
linear} if there is a finite partition $\Pm$ of $X$ such that all
restrictions $f|_A$, $A\in \Pm$ are linear. We speak analogously
of linear and piecewise linear maps $g:X\to G^n$.
 \end{definition}
The following definition fixes the notion of (Presburger) cells.
 \begin{definition}\label{def cell}
A cell of type $(0)$ (also called a $(0)$-cell) is a point
$\{a\}\subset G$. A $(1)$-cell is a set with infinite cardinality
of the form
 \begin{equation}\label{cell dim 1}
  \{x\in G\mid \alpha\sq_1 x\sq_2 \beta,\ x\equiv c \pmod{n}\},
 \end{equation}
with $\alpha,\beta$ in $G$, integers $0\leq c<n$ and $\sq_i$
either $\leq$ or no condition. Let $i_j\in\{0,1\}$ for
$j=1,\ldots,m$ and $x=(x_1,\ldots,x_m)$. A
$(i_1,\ldots,i_{m},1)$-cell is a set $A$ of the form
 \begin{equation}\label{cell}
  A = \{(x,t)\in G^{m+1} \mid x\in D,\ \alpha(x)\sq_1 t \sq_2
  \beta(x),\ t\equiv c \pmod{n}\},
 \end{equation}
with $D=\pi_{m}(A)$ a $(i_1,\ldots,i_m)$-cell,
 $\alpha,\ \beta:D\to G$ linear functions, $\sq_i$ either
$\leq$ or no condition and integers $0\leq c<n$ such that the
cardinality of the fibers $A_x=\{t\in G\mid (x,t)\in A\}$ can not
be bounded uniformly in $x\in D$ by an integer.
 \\
 A $(i_1,\ldots,i_{m},0)$-cell is a set of the form
 \[
\{(x,t)\in G^{m+1} \mid x\in D,\ \alpha(x)=t\},
 \]
with $\alpha:D\to G$ a linear function and $D\subset G^{m}$ a
$(i_1,\ldots,i_m)$-cell.
 \end{definition}
 \begin{remarks}\label{remark cell}
\item[(i)]
Although we consider in Definition \ref{def cell} a condition on
the cardinality of fibers, the type of a cell does not alter if
one takes
elementary extensions.\\
\item[(ii)]\label{remark cell item}
 To an infinite $(i_1,\ldots,i_m)$-cell $A\subset G^m$ we can associate
(as in \cite{vdD}) a projection $\pi_A:G^m\to G^k$ such that the
restriction of  $\pi_A$ to $A$ gives a bijection from $A$ onto a
$(1,\ldots,1)$-cell $A'\subset G^k$. Also,  a $(i_1, \ldots,
i_m)$-cell
is finite if and only if $i_1 = \cdots = i_m = 0$, and then it is a singleton.\\
\item[(iii)]
Let $A$ be a $(i_1,\ldots,i_{m},1)$-cell as in Eq.~(\ref{cell}),
then it is clear that a linear function $f:A\to G$ can be written
as
 \begin{equation}\label{lin dim 1}
 f(x,t)=a(\frac{t-c}{n})+\gamma(x),\quad  (x,t)\in A,
 \end{equation}
with $a$ an integer, $\gamma:D\to G$ a linear function and $c,n,D$
as in Eq.~(\ref{cell}).
 \end{remarks}
 \begin{theorem}[Cell Decomposition]\label{cell decomp}
Let $X\subset G^m$ and $f:X\to G$ be $\Lp$-definable. Then there
exists a finite partition $\Pm$ of $X$ into cells, such that the
restriction $f|_A:A\to G$ is linear for each cell $A\in\Pm$.
Moreover, if $X$ and $f$ are $S$-definable, then also the parts
$A$ can be taken $S$-definable.
 \end{theorem}
 \begin{proof}[Proof by induction on $m$.]
If $X\subset G$, $f:X\to G$ are $\Lp$-definable, then Theorem
\ref{cell decomp} follows easily by using quantifier elimination
and elementary properties of linear congruences. Alternatively,
the more general Thm.~4.8 of \cite{Point} can be used to prove
this one dimensional version (see also Proposition~\ref{piecewise
linear} below).
 Let $X\subset G^{m+1}$ and  $f:X\to G$ be
$\Lp$-definable, $m>0$. We write $(\sq_1, \sq_2)\in
\{\leq,\emptyset\}^2$ to say that $\sq_1$, resp.~$\sq_2$,
represents either the symbol $\leq$ or no condition. Let  $\SSS$
be the set $\Z\times\{(n,c)\in\Z^2 \mid 0\leq c <
n\}\times\{\leq,\emptyset\}^2$.
 For any $d=(a,n,c,\sq_1, \sq_2) \in \SSS$
and $\xi=(\xi_1,\xi_2,\xi_3) \in G^3$ we define a Presburger
function
\[
 F_{(d,\xi)}:\{t\in G\mid \xi_1 \sq_1 t \sq_2 \xi_2,\ t\equiv c
\pmod{n}\}\to G:t\mapsto a(\frac{t-c}{n})+\xi_3.
\]
The domain  $Dom(F_{(d,\xi)})$ of such a function $F_{(d,\xi)}$ is
either empty, a $(1)$-cell or a finite union of $(0)$-cells. For
fixed $k>0$ and  $d\in \SSS^k$, let $\varphi_{(d,k)}(x,\xi)$ be a
Presburger formula in the free variables $x=(x_1,\ldots,x_m)$ and
$\xi=(\xi_1,\ldots,\xi_k)$, with
$\xi_i=(\xi_{i1},\xi_{i2},\xi_{i3})$, such that $G\models
\varphi_{(d,k)}(x,\xi)$ if and only if the following are true:
 \begin{itemize}
\item[(i)] $x\in\pi_m(X)$,\\
\item[(ii)] the collection  of the domains $Dom(F_{(d_i,\xi_i)})$ for $i=1,\ldots,k$ forms
 a partition of the fiber $X_x\subset
G$,\\
\item[(iii)] $F_{(d_i,\xi_i)}(t)=f(x,t)$ for
each  $t\in Dom(F_{(d_i,\xi_i)})$ and $i=1,\ldots,k$.
 \end{itemize}
Now we define for each $k$ and $d \in \SSS^k$ the set
\[
B_{(d,k)}= \{x\in G^m\mid \exists \xi\quad
\varphi_{(d,k)}(x,\xi)\}.
\]
Each set $B_{(d,k)}$ is $\Lp$-definable and the (countable)
collection $\{B_{(d,k)}\}_{k,d}$ covers $\pi_m(X)$ since each
$x\in\pi_m(X)$ is in some $B_{(d,k)}$ by the induction basis. We
can do this construction in any elementary extension of $G$, so by
logical compactness we must have that finitely many sets of the
form $B_{(d,k)}$ already cover $G^m$. Consequently, we can take
Presburger sets $D_1,\ldots,D_s$ such that $\{D_i\}$ forms a
partition of $G^m$ and each $D_i$ is contained in a set
$B_{(d,k)}$ for some $k$ and $k$-tuple $d$. For each
$i=1,\ldots,s$, fix a $k$ and $k$-tuple $d$ with $D_i\subset
B_{(d,k)}$, then we can define the Presburger set
\[
\Gamma_i=\{(x,\xi)\in D_i\times G^{3k}\mid
\varphi_{(d,k)}(x,\xi)\}
\]
satisfying $\pi_m(\Gamma_i)=D_i$ by construction. Since the theory
$\mathrm{Th}(G,\Lp)$ has definable Skolem functions, we can choose
definably for each $x\in D_i$ tuples $\xi\in G^{3k}$ such that
$(x,\xi)\in\Gamma_i$. Combining it all, it follows that there
exists a finite partition $\Pm$ of $X$ consisting of Presburger
sets of the form
\[
A=\{(x,t)\in G^{m+1}\mid x\in C,\ \alpha(x)\sq_1 t \sq_2
\beta(x),\ t\equiv c \pmod{n}\},
\]
such that $f|_A$ maps $(x,t)\in A$ to
$a(\frac{t-c}{n})+\gamma(x)$; with $\alpha,\beta,\gamma:C\to G$
and $C\subset G^m$ $\Lp$-definable, $\sq_i$ either $\leq$ or no
condition, integers $a$, $0\leq c<n$ and $\pi_m(A)=C$. The theorem
now follows after applying the induction hypothesis to $C$ and
$\alpha,\beta,\gamma:C\to G$ and partitioning further.
\end{proof}
\section{Dimension theory for Presburger arithmetic}
Any Presburger structure $(G,\Lp)$ satisfies the exchange property
for algebraic closure. This is a corollary of a more general
result in \cite{Wagner} but can also be proven  using the cell
decomposition theorem elementarily. In particular this yields an
algebraic dimension
function on the Presburger sets in the following (standard) way.\\
 \begin{definition}
Let $X\subset G^m$ be $A$-definable for some finite set $A$ by a
formula $\varphi(x,a)$ where $a=(a_1,\ldots,a_s)$ enumerates $A$,
then the (algebraic) dimension of $X$, written $\mathrm{dim}(X)$,
is the greatest integer $k$ such that in some elementary extension
$\bar G$ of $G$ there exists $x=(x_1,\ldots,x_m)\in \bar G^{m}$
with $\bar G\models \varphi(x,a)$ and
$\mathrm{rk}(x_1,\ldots,x_m,a_1,\ldots,a_s)-\mathrm{rk}(a_1,\ldots,a_s)=k$,
where $\mathrm{rk}(B)$ of a set $B\subset \bar G$ is the
cardinality of a maximal algebraically independent subset of $B$
(in the sense of model theory, see \cite{Hodges}).
\end{definition}
This dimension function is independent of the choice of a set of
defining parameters $A$ and  the following properties of algebraic
dimension are standard.
 \begin{proposition}\label{behaves well}
\item[(i)]For Presburger sets $X,Y\subset G^m$ we have
$\mathrm{dim}(X\cup Y)=\max(\mathrm{dim}X,\mathrm{dim}Y)$.\\
\item[(ii)]Let $f:X\to G^m$ be $\Lp$-definable, then $\mathrm{dim}(X)\geq
\mathrm{dim}(f(X))$.
 \end{proposition}
The dimension of a cell $C$ is directly related to the type of $C$
(see Lemma \ref{algdim cell}). Also, if we have a Presburger set
$X$ and a finite partition $\Pm$ of $X$ into cells, the dimension
of $X$ is directly related to the types of the cells in $\Pm$ (see
Cor.~\ref{dim, partitie}).
 \begin{lemma}\label{algdim cell}
Let $C\subset G^m$ be a $(i_1,\ldots,i_m)$-cell, then
$\mathrm{dim}(C)=i_1+\ldots+i_m$.
 \end{lemma}
 \begin{proof} For a $(0)$- and a $(1)$-cell this is clear.
Possibly after projecting, we may suppose that $C\subset G^m$ is a
$(1,\ldots,1)$-cell. The Lemma follows now from the definition of
the type of a cell using induction on $m$ and a compactness
argument.
 \end{proof}
 \begin{cor}\label{dim, partitie}
For any Presburger set $X\subset G^m$ and any finite partition
$\Pm$ of $X$ into cells we have
 \begin{eqnarray}
 \mathrm{dim}(X)
  & = &
  \max\{i_1+\ldots+i_m\mid C\in\Pm,\mbox{ $C$ is a
  $(i_1,\ldots,i_m)$-cell}\}\nonumber
  \\
  & = &
  \max\{i_1+\ldots+i_m\mid X \mbox{ contains a
  $(i_1,\ldots,i_m)$-cell}\}.\label{Eqdim1}
 \end{eqnarray}
 \end{cor}
 \begin{proof} The first equality is a consequence of Lemma~\ref{algdim cell}
and Proposition~\ref{behaves well}. To prove (\ref{Eqdim1}) we
take a $(i_1,\ldots,i_m)$-cell $C\subset X$ such that
$i_1+\ldots+i_m$ is maximal. By the cell decomposition we can
obtain a partition $\Pm$ of $X$ into cells such that $C\in\Pm$.
Now use the previous equality to finish the proof.
 \end{proof}
 \begin{remark}
It is also possible to take Eq.~(\ref{Eqdim1}) as the definition
for the dimension of a Presburger set and to proceed similarly as
in \cite{vdD} by van den Dries to develop a dimension theory for
Presburger structures.
 \end{remark}
 \section{Classification of Presburger sets}\label{sectionclassification}
The cell decomposition theorem provides us with the technical
tools to classify the $\emptyset$-definable Presburger sets up to
$\Lp$-definable bijection. The key step to this classification is
a rectilinearisation theorem, which also has a parametric
formulation. We recall that $G$ denotes a $Z$-group and $H=\{x\in
G\mid x\geq0\}$, we write $H^0=\{0\}$. Also notice that a set $A$
is $\emptyset$-definable if and only if $A$ is $\Z$-definable, to
be precise, definable over $\Z\cdot 1\subset G$.
 \begin{theorem}[Rectilinearisation]\label{recti}
Let $X$ be a $\emptyset$-definable Presburger set, then there
exists a finite partition $\Pm$ of $X$ into $\emptyset$-definable
Presburger sets, such that for each $A\in \Pm$ there is an integer
$l\geq 0$ and a $\emptyset$-definable linear bijection $f:A\to
H^l$.
 \end{theorem}
 \begin{proof} We give a proof by induction on $\mathrm{dim}X$. If $\mathrm{dim} X=0$ then
$X$ is finite and the theorem follows, so we choose a Presburger
set $X$ with $\mathrm{dim} X=m+1$, $m\geq0$. Any $\Lp$-definable
object occurring in this proof will be $\emptyset$-definable; we
will alternately apply $\emptyset$-definable linear bijections and
partition further. By the cell decomposition theorem and possibly
after projecting (see the remark after Definition \ref{def cell}),
we may suppose that $X$ is a $(1,\ldots,1)$-cell contained in
$G^{m+1}$, so we can write
\[
X=\{(x,t)\in G^{m+1}\mid x\in D,\ \alpha(x)\sq_1 t\sq_2 \beta(x),\
t\equiv c\pmod{n}\},
\]
with $x=(x_1,\ldots,x_m)$, $\pi_m(X)=D\subset G^m$ a
$(1,\ldots,1)$-cell, integers $0\leq c<n$, $\alpha,\beta:D\to G$
$\emptyset$-definable linear functions and $\sq_i$ either $\leq$
or no condition. By induction we may suppose that $D=H^m$. If both
$\sq_1$ and $\sq_2$ are no condition, the theorem follows easily,
so we may suppose that one of the $\sq_i$, say $\sq_1$, is $\leq$.
Moreover, after a linear transformation
$(x,t)\mapsto(x,\frac{t-c}{n})$ we may assume that $c=0$ and
$n=1$, then we can apply the following linear bijection
\[
   f:X \to A:(x,t) \mapsto (x_1,\ldots,x_m,t-\alpha(x)),
\]
onto
 \begin{equation*}
A  =  \{(x,t)\in H^{m+1}\mid t \sq_2 \beta(x)-\alpha(x)\}.
 \end{equation*}
Because $\beta(x)-\alpha(x)$ is a linear function from $H^m$ to
$G$ there are integers $k_i$ such that
\begin{equation}\label{some cell}
A  =  \{(x,t)\in H^{m+1}\mid t \sq_2 k_0+\sum_{i=1}^m k_i x_i\}.
\end{equation}
Moreover, since $\pi_m(A)=H^m$, all integers $k_i$ must be
nonnegative. We proceed by induction on $k_1\geq0$. If $k_1=0$
then $A=H\times\{(x_2,\ldots,x_m,t)\in H^m\mid t\leq
k_0+\sum_{i=2}^m k_i x_i\}$ and the theorem follows by induction
on the dimension. Now suppose $k_1>0$, then we partition $A$ into
two parts
 \begin{eqnarray*}
 A_1 & = & \{(x,t)\in H^{m+1}\mid t\leq x_1-1\},\\
 A_2 & = & \{(x,t)\in H^{m+1}\mid x_1\leq t\leq k_0+\sum_{i=1}^m k_i x_i\},
 \end{eqnarray*}
where $\pi_m(A_2)= H^m$ and $\pi_m(A_1)=\{x\in H^m \mid 1\leq
x_1\}$. We apply the linear bijection
\[
A_2\to B:(x,t)\mapsto (x_1,\ldots,x_m,t-x_1)
\]
with
\[
B=\{(x,t)\in H^{m+1}\mid  t\leq k_0+(k_1-1)x_1+\sum_{i=2}^m k_i
x_i\}
\]
and the theorem  for $B$ follows by induction on $k_1$. We
conclude the proof by the following linear bijection:
\[
A_1\to H^{m+1}:(x,t)\mapsto(x_1-1-t,x_2,\ldots,x_m,t).
\]
 \end{proof}
 \begin{theorem}[Parametric Rectilinearisation]\label{param recti}
Let $X\subset G^{m+n}$ be a $\emptyset$-definable Presburger set,
then there exists a finite partition $\Pm$ of $X$ into
$\emptyset$-definable Presburger sets, such that for each $A\in
\Pm$ there is a set $B\subset G^{m+n}$ with $\pi_m(A)=\pi_m(B)$
and a $\emptyset$-definable family $\{f_\lambda\}_{\lambda\in
\pi_m(A)}$ of linear bijections $f_\lambda:A_\lambda\subset G^n\to
B_\lambda\subset G^n$ with $B_\lambda$ a set of the form
$H^l\times\Lambda_\lambda$ where $\Lambda_\lambda$ is a bounded
$\lambda$-definable set and the integer $l$ only depends on
$A\in\Pm$.
 \end{theorem}
 \begin{proof} We give a proof by induction on $n$, following the lines
of the proof of Theorem \ref{recti}. So we assume that $X$ is a
cell
\[
X=\{(\lambda,x,t)\in G^{m+(n+1)}\mid (\lambda,x)\in D,\
\alpha(\lambda,x)\sq_1 t\sq_2 \beta(\lambda,x),\ t\equiv
c\pmod{n}\},
\]
with $\lambda=(\lambda_1,\ldots,\lambda_m)$, $x=(x_1,\ldots,x_n)$,
$D\subset G^{m+n}$ a cell, integers $0\leq c<n$,
$\alpha,\beta:D\to G$ $\emptyset$-definable linear functions and
$\sq_i$ either $\leq$ or no condition. By subsequently applying
the induction hypothesis to $D$, partitioning further and applying
linear bijections (similar as to obtain Eq.~(\ref{some cell}) in
the proof of Theorem \ref{recti}, keeping the parameters $\lambda$
fixed now), we may assume that $X$ has the form
 \[
X=\{(\lambda,x,t)\in G^{m+n+1}\mid (\lambda,x)\in D',\ 0\leq t
\leq \gamma(\lambda,x)\},
 \]
with $\pi_{m+n}(X)=D'\subset G^{m+n}$ a Presburger set such that
for each $\lambda\in\pi_m(D')$ $D'_\lambda=H^l\times
\Gamma_\lambda$ where $\Gamma_\lambda$ is a $\lambda$-definable
bounded set, $l$ a fixed positive integer and $\gamma:D'\to G$ a
$\emptyset$-definable linear function. If $l=0$, $X_\lambda$ is a
bounded set for each $\lambda$ and the theorem follows
immediately. Let thus $l\geq 1$, i.e.~the projection of $X$ on the
$x_1$-coordinate is $H$, then the function $\gamma$ can be written
as $(\lambda,x)\mapsto k_1x_1+\gamma'(\lambda,x_2,\ldots, x_m)$
with $k_1$ an integer, necessarily nonnegative because the
projection of $X$ on the $x_1$-coordinate is $H$ and $\gamma'$ is
a linear function. The reader can finish the proof by induction on
$k_1\geq 0$, similar as in the proof of Theorem \ref{recti}.
 \end{proof}
 \begin{theorem}[Classification]\label{classification}
Let $X$ be a $\emptyset$-definable Presburger set with
$\mathrm{dim}X=m>0$, then there exists a $\emptyset$-definable
Presburger bijection $f:X\to G^m$. In other words, there exists a
$\emptyset$-definable Presburger bijection between two infinite
$\emptyset$-definable Presburger sets $X,Y$ if and only if
$\mathrm{dim} X=\dim Y$.
 \end{theorem}
 \begin{proof} Let $X$ be $\emptyset$-definable and infinite.  We use induction on
$\mathrm{dim}X=m$. We say for short that two Presburger sets $X,Y$
are isomorphic if there exists a $\emptyset$-definable Presburger
bijection between them and write $X\cong Y$. If $m=1$, then
Theorem \ref{recti} yields a partition $\Pm$ of $X$ such that each
part is either a point or isomorphic to $H$. Consider the
bijections
\[
 \begin{array}{l}
f_1:H\to G :\left\{\begin{array}{rcl} 2x & \mapsto & x,
 \\ 2x+1 & \mapsto& -x,\end{array}\right.\\

f_2:H\cup\{-1\}\to H:x\mapsto x+1,\\

f_3:(\{0\}\times H) \cup (\{1\}\times H)\to
H:\left\{\begin{array}{rcl} (0,x) & \mapsto & 2x,\\ (1,x) &
\mapsto& 2x+1;\end{array}\right.
 \end{array}\]
the bijections $f_1,f_2$, applied repeatedly to (isomorphic copies
of) parts in $\Pm$ yield a definable bijection from $X$ onto $H$
and thus $G\cong X$ by applying $f_1$ (in the obvious way). Now
let $\mathrm{dim}X=m>1$. Using Theorem \ref{recti} we find a
partition $\Pm$ of $X$ such that each part is isomorphic to $H^l$
and thus to $G^l$ since $H\cong G$ by $f_1$. Since
$\mathrm{dim}X=m$, at least one part is isomorphic to $G^m$. Take
$A,B\in\Pm$ with $A\cong G^m$ and $B\cong G^l$, then it suffices
to prove that $A\cup B\cong G^m$. If $l=0$ this is clear and if
$l>0$ then $A\cup B\cong G\times(A'\cup B')$ for some disjoint and
$\emptyset$-definable sets $A',B'$ with $A'\cong G^{m-1}$ and
$B'\cong G^{l-1}$. The induction hypothesis applied to $A'\cup B'$
finishes the proof.
 \end{proof}
 \section{Expansions of $Z$-groups}\label{section Michaux}
We define the notion of Presburger minimality ($\Lp$-minimality)
for expansions of Presburger structures $(G,\Lp)$. This notion of
$\Lp$-minimality is a concrete instance of the general notion of
$\Lm$-minimality as in \cite{Macpherson} and has already been
studied in \cite{Point}.
 \begin{definition}
Let $G$ be a $Z$-group and $\Lm$ an expansion of the language
$\Lp$, then we say that $(G,\Lm)$ is $\Lp$-minimal if every
$\Lm$-definable subset of $G$ is already  $\Lp$-definable
(allowing parameters as always). We say that $\mathrm{Th}(G,\Lm)$
is $\Lp$-minimal if every model of this theory is $\Lp$-minimal.
 \end{definition}
Comparing this notion with the terminology of \cite{Point}, a
structure $(G,\Lm)$ is $\Lp$-minimal if and only if it is a
\emph{discrete coset-minimal group without definable proper convex
subgroups} (see \cite{Point}). Theorem 4.8 of \cite{Point} says
that a definable function in one variable between such groups is
piecewise linear. We reformulate this result with our terminology.
 \begin{proposition}[\cite{Point}, Thm.~4.8]\label{piecewise linear}
Let $(G,\Lm)$ be $\Lp$-minimal, then any definable function
$f:G\to G$ is piecewise linear.
 \end{proposition}
Proposition \ref{piecewise linear} allows us to repeat without any
change the compactness argument of the proof of the cell
decomposition theorem for any model of a $\Lp$-minimal theory.
This leads to the following description of $\Lp$-minimal theories.
 \begin{theorem}\label{strong}
Let $(G,\Lm)$ be an expansion of a Presburger structure $(G,\Lp)$,
then the following are equivalent:
 \item[(i)] $\mathrm{Th}(G,\Lm)$ is $\Lp$-minimal;
 \item[(ii)] $(G,\Lm)$ is a definitional expansion of $(G,\Lp)$; precisely,
 any $\Lm$-definable set $X\subset G^m$ is already $\Lp$-definable.\\
Thus, the theory $\mathrm{Th}(G,\Lp)$ does not admit any proper
$\Lp$-minimal expansion.
 \end{theorem}
 \begin{proof} Any Presburger minimal theory has definable Skolem
functions. For if $X\subset G^{m+1}$ is a definable set in some
model $G$, we can choose definably for any $x\in\pi_m(X)$  the
smallest nonnegative element in $X_x$ if there is any, and the
largest  negative element otherwise. This implies the definability
of Skolem functions by induction. Now replace in the statement of
the cell decomposition Theorem (theorem \ref{cell decomp}) the
word \emph{$\Lp$-definable} by \emph{$\Lm$-definable}. Then repeat
the case $m=1$ of the proof of Theorem \ref{cell decomp}, using
now the $\Lp$-minimality and Proposition \ref{piecewise linear}.
Using the same compactness argument as in the proof of Theorem
\ref{cell decomp} we find that any $\Lm$-definable set $X\subset
G^m$ is a finite union of Presburger cells, thus a fortiori, $X$
is $\Lp$-definable.
 \end{proof}
 \begin{remark}
For an arbitrary expansion $(G,\Lm)$ of $(G,\Lp)$ it is, as far as
I know, an open problem whether the statements (i) and (ii) of
Thm.~\ref{strong} are equivalent with the following:
\item[(iii)] $(G,\Lm)$ is $\Lp$-minimal.\\
In the special case  $G=\Z$, statements (i), (ii) and (iii) are
indeed equivalent, proven by C.~Michaux and R.~Villemaire
in~\cite{Michaux}.
 \end{remark}
 \section{Application to p-minimal fields}\label{section p}
In this section, we let $K$ be a $p$-adically closed field with
value group $G$. Recall that a $p$-adically closed field is a
field $K$ which is elementary equivalent to a finite field
extension of the field $\Q_p$ of $p$-adic numbers; in particular,
the value group $G$ is a $Z$-group and $K$ has quantifier
elimination in the Macintyre language
$\Lmac=\langle+,-,.,0,1,\{P_n\}_{n\geq 1}\rangle$ where $P_n$
denotes the set of $n$-th powers in $K^\times$. We write $v:K\to
G\cup\{\infty\}$ for the valuation map and for any $m>0$ we write
$\bar v$ for the map $\bar v:(K^\times)^ m\to G^m:x\mapsto
(v(x_1),\ldots,v(x_m))$. We give a definition of p-minimality,
extending the original definition of \cite{Haskell} slightly.
 \begin{definition}
Let $K$ be a $p$-adically closed field and let $(K,\Lm)$ be an
expansion of $(K,\Lmac)$. We say that the structure $(K,\Lm)$ is
p-minimal if any $\Lm$-definable subset of $K$ is already
$\Lmac$-definable (allowing parameters). The theory
$\mathrm{Th}(K,\Lm)$ is called p-minimal if every model of this
theory is p-minimal.
 \end{definition}
Examples of p-minimal fields known at this moment are $p$-adically
closed fields with the semi-algebraic structure and with
subanalytic structure with restricted power series (see
\cite{vdDHM}). Theorem \ref{application} exhibits a close
connection between definable sets in a p-minimal field $K$ and
Presburger sets in the value group $G$ of $K$; to prove it, we use
Lemma \ref{Mac}, which is a reformulation of the interpretability
of $(G,\Lp)$ in $(K,\Lmac)$.
 \begin{lemma}\label{Mac}
Let $K$ be a $p$-adically closed field with value group $G$, then
for any $\Lp$-definable set $S\subset G^m$ the set $\bar
v^{-1}(S)=\{(x_1,\ldots,x_m)\in (K^\times)^m\mid\bar v(x)\in S\}$
is $\Lmac$-definable.
 \end{lemma}
 \begin{proof} Let $S\subset G^m$ be $\Lp$-definable. By Theorem \ref{cell
decomp} we may suppose that $S$ is a Presburger cell. The Lemma
follows now inductively  from the fact that conditions imposed on
$(x_1,\ldots,x_{m-1},t)\in (K^\times)^m$ of the form $\pm
v(t)\leq\frac{1}{e}(\sum_{i=1}^{m-1}a_iv(x_i))+d$ or $v(t)\equiv
c\pmod{n}$ are $\Lmac$-definable for any integers $a_i,e$,
$e\not=0$, $0\leq c< n$ and $d\in G$ (see e.g.~\cite[Lemma
2.1]{Denef2}).
 \end{proof}
 \begin{theorem}\label{application}
Let $(K,\Lm)$ be a p-minimal field with p-minimal theory and let
$G$ be the value group of $K$. Then for any $\Lm$-definable set
$X\subset (K^\times)^m$ the set
\[
\bar v(X)=\{(v(x_1),\ldots,v(x_m))\in G^m\mid (x_1,\ldots,x_m)\in
X\} \subset G^m
\]
is $\Lp$-definable.
 \end{theorem}
 \begin{proof} Put $S_m=\{\bar v(X)\subset G^m\mid X\subset (K^\times)^m, \mbox{ $X$
is $\Lm$-definable}\}$, then it is easy to see that the collection
$(S_m)_{m\geq 0}$  determines  a  structure on $G$ (i.e.~the
collection $\cup_m S_m$ is precisely the collection of
$\Lm'$-definable sets for some language $\Lm'$). We first argument
that this structure is in fact $\Lp$-minimal. Choose a
$\Lm$-definable set $X\subset K^\times$, then, by p-minimality,
$X$ is $\Lmac$-definable. We can thus apply the $p$-adic
semi-algebraic cell decomposition (\cite{Denef2}, in the
formulation of \cite[Lemma 4]{C})
 to the set $X$ to obtain that $X$ is a
finite union of $p$-adic cells, i.e.~sets of the form
\[
\{x\in K\mid  v(a_1)\sq_1 v(x-c)\sq_2 v(a_2),\ x-c\in\lambda
P_n\}\subset K^\times,
\]
with $a_1,a_2,c,\lambda\in K$ and $\sq_i$ either $\leq,<$ or no
condition. The image under $v$ of such a cell is either a finite
union of $(0$)-cells or a $(1)$-cell and thus a $\Lp$-definable
subset of $G$. By consequence,  the structure $(S_m)_{m\geq_0}$ is
$\Lp$-minimal. By the Presburger minimality of $(S_m)_{m\geq 0}$,
the p-minimality of $\mathrm{Th}(K,\Lm)$ and Lemma \ref{Mac} to
interprete $G$ into $K$, we can repeat the compactness argument of
the proof of the cell decomposition theorem \ref{cell decomp} for
the structure $(S_m)_m$ on $G$ to find that each $A\in \cup_m S_m$
is a finite union of Presburger cells. This proves the theorem.
 \end{proof}
 \begin{remark}
For a $p$-adically closed field $K$ it is proven by the author in
\cite{C} that there exists a $\Lmac$-definable bijection $X\to Y$
between two infinite parameter free $\Lmac$-definable sets $X,Y$
if and only if $\dim X=\dim Y$. This is proven by reducing to a
Presburger problem similar to the classification theorem in
section \ref{sectionclassification}. An analogous classification
for $\emptyset$-definable sets in arbitrary p-minimal structures
is not known up to now.
 \end{remark}
 \section{Elimination of imaginaries}
As a last application of the cell decomposition theorem we prove
uniform elimination of imaginaries for Presburger structures. We
say that a structure $(M,\Lm)$ has uniform elimination of
imaginaries if for any $\emptyset$-definable equivalence relation
on $M^k$ there exists a $\emptyset$-definable function $F:M^k\to
M^r$ for some $r$ such that two tuples $x,y\in M^k$ are equivalent
if and only if $F(x)=F(y)$.
 \begin{theorem}
 The theory $Th(\Z,\Lp)$ has uniform elimination of imaginaries,
precisely, any Presburger structure $(G,\Lp)$ eliminates
imaginaries uniformly.
 \end{theorem}
 \begin{proof} Since $\mathrm{Th}(\Z,\Lp)$ has definable Skolem functions,
we only have to prove the following statement for an arbitrary
$Z$-group $G$ (see e.g.~\cite[Lemma 4.4.3]{Hodges}). For any
$\emptyset$-definable Presburger set $X\subset G^{m+1}$ there
exists a $\emptyset$-definable Presburger function $F:G^m\to G^n$
for some $n$, such that $F(x)=F(x')$ if and only if $X_x=X_{x'}$
(if $x\not\in \pi_m(X)$ then we put $X_x=\emptyset$). So let
$X\subset G^{m+1}$ be a $\emptyset$-definable Presburger set.
Apply the cell decomposition theorem to obtain a partition $\Pm$
of $X$ into cells. For each cell $A\in \Pm$ of the form
$A=\{(x,t)\in G^{m+1}\mid x\in D,\ \alpha(x)\sq_{1A} t\sq_{2A}
\beta(x),\ t\equiv c\pmod{n}\}$ (as in Eq.~\ref{cell}) and each
$\xi=(\xi_1,\xi_2)\in G^2$ we define a set
\[
C_A(\xi)=\{t\in G\mid \xi_1\sq_{1A} t\sq_{2A} \xi_2,\ t\equiv
c\pmod{n}\}.
\]
Notice that for each $x\in \pi_m(X)$ we have at least one
partition of $X_x$ into sets of the form $C_A(\xi)$ with $A\in\Pm$
and $\xi\in G^2$. For $x,y\in G$ we write $x\lhd y$ if and only if
one of the following conditions is satisfied
 \begin{itemize}
\item[(i)] $0\leq x< y$,\\
\item[(ii)]$0< x\leq -y$,\\
\item[(iii)]$0<-x<y$,\\
\item[(iv)]$0<-x<-y$.
 \end{itemize}
This gives a new ordering $0 \lhd 1\lhd -1\lhd 2\lhd -2\lhd\ldots$
on $G$  with zero as its smallest element. For each $k>0$ we also
write $\lhd$ for the lexicographical order on $G^k$ built up with
$\lhd$. The order $\lhd$ is $\Lp$-definable and each Presburger
set has a unique $\lhd$-smallest element. For each $x\in G^m$ and
each $I\subset \Pm$ with cardinality $|I|=s\geq0$ we let
$y_I(x)=(\xi_A)_{A\in I}$, $\xi_A=\in G^2$, be the $\lhd$-smallest
tuple in $G^{2s}$ such that $\cup_{A\in I}C_A(\xi_A)=X_x$ if there
exists at least one such tuple and we put $y_I(x)=(0,\ldots,0)\in
G^{2s}$ otherwise. One can reconstruct the set $X_x$ given all
tuples $y_I(x)$, $I\subset \Pm$. Let $F$ be the function mapping
$x\in \pi_m(X)$ to $y=(y_I(x))_{I\subset \Pm}$. Since the
lexicographical order $\lhd$ is $\Lp$-definable it is clear that
$F$ is $\Lp$-definable and that $F(x)=F(x')$ if and only if
$X_x=X_{x'}$ for each $x,x'\in G^m$.
 \end{proof}
\subsection*{Acknowledgment}
I would like to thank J.~Denef, D.~Haskell, C.~Michaux, F.~Point, F.~
Wagner and K.~Zahidi for fruitful discussions during the preparation of
this paper.
\bibliography{presbib}
\bibliographystyle
{asl}
\end{document}